\documentclass[11pt]{article}
\usepackage[utf8]{inputenc}
\usepackage[a4paper]{geometry}
\usepackage[english]{babel}
\usepackage[super]{nth}
\usepackage{amsmath, amsfonts, amsthm, physics, enumitem,amssymb}
\usepackage{authblk}
\usepackage{titlesec}
\usepackage{tikz} 
\usetikzlibrary{positioning}
\usepackage{mathtools}

\usetikzlibrary{patterns,arrows,decorations.markings}

\usepackage{tikz-cd}
\usepackage{centernot}

\makeatletter
\newtheorem*{rep@theorem}{\rep@title}
\newcommand{\newreptheorem}[2]{%
\newenvironment{rep#1}[1]{%
 \def\rep@title{#2 \ref{##1}}%
 \begin{rep@theorem}}%
 {\end{rep@theorem}}}
\makeatother

\newtheorem{thm}{Theorem}[section]
\newreptheorem{thm}{Theorem}

\newtheorem{cor}[thm]{Corollary}

\newtheorem{q}[thm]{Question}
\newreptheorem{q}{Question}
\newtheorem{prop}[thm]{Proposition}
\newtheorem{defn}[thm]{Definition}
\newtheorem{lem}[thm]{Lemma}

\theoremstyle{remark}

\theoremstyle{definition}
\newtheorem{ex}[thm]{Example}

\newcommand{\der}[1]{\frac{\partial}{\partial#1}}
\newcommand{\dif}[1]{\mathrm{d}#1}
\newcommand{\Nil}{\mathrm{Nil}}
\newcommand{\orb}{\mathrm{Orb}}
\newcommand{\im}{\mathrm{im} \,} 

\makeatletter
\def\thanks#1{\protected@xdef\@thanks{\@thanks
        \protect\footnotetext{#1}}}
\makeatother

\title{Bott-Chern and $\bar{\partial}$ Harmonic Forms on Almost Hermitian 4-Manifolds}
\author{Tom Holt \thanks{Thomas.Holt@warwick.ac.uk}}
\affil{Mathematics Institute,  University of Warwick, Coventry, CV4 7AL, UK}
\date{}

\begin{document}

\maketitle

\begin{abstract}
    We prove that on a compact almost Hermitian 4-manifold the space of $\bar\partial$-harmonic $(1,1)$-forms always has dimension $h_{\bar\partial}^{1,1} = b_- +1$ or $b_-$, whilst the space of Bott-Chern harmonic $(1,1)$-forms always has dimension $h_{BC}^{1,1} = b_- +1$. We also perform calculations of $h^{2,1}_{BC}$ and $h^{1,2}_{BC}$ on the Kodaira-Thurston manifold, thereby providing a full account of when $h^{p,q}_{BC}$ is or is not invariant of the choice of almost Hermitian metric. Finally, we introduce a decomposition of the space of $L^2$ functions on all torus bundles over $S^1$, which has proven useful for solving linear PDEs, and we demonstrate its use in the calculation of $h^{p,q}_{\bar\partial}$.  
\end{abstract}
\tableofcontents

\section{Introduction}

On an almost complex manifold $(X,J)$ endowed with an almost Hermitian metric, we can define the spaces 
$$\mathcal{H}^{p,q}_{\bar\partial} := \ker \Delta_{\bar\partial}|_{\mathcal{A}^{p,q}}\quad \quad \mathcal{H}^{p,q}_{BC} := \ker \Delta_{BC}|_{\mathcal{A}^{p,q}} $$
of $\bar\partial$-harmonic and Bott-Chern harmonic $(p,q)$-forms respectively. The corresponding laplacians are both elliptic differential operators and so the dimension of these spaces, denoted by $h^{p,q}_{\bar\partial}$ and $h^{p,q}_{BC}$, is finite whenever $X$ is compact.

In the case when $J$ is integrable, \textit{i.e.} $J$ arises from a complex structure on $X$, $\mathcal{H}^{p,q}_{\bar\partial}$ is isomorphic to the Dolbeault cohomology
$$H^{p,q}_{\bar\partial}:= \frac{\ker \bar\partial}{\im \bar\partial} $$
whilst $\mathcal{H}_{BC}^{p,q}$ is isomorphic to the Bott-Chern cohomology
$$H^{p,q}_{BC}:= \frac{\ker \partial \cap \ker \bar\partial}{\im \partial\bar\partial }. $$ 
One consequence of this is that on complex manifolds $h^{p,q}_{\bar\partial}$ and $h^{p,q}_{BC}$ are both independent of the choice of Hermitian metric. However, when $J$ is non-integrable it is no longer the case that $\bar{\partial}^2 = 0 $ and so the cohomology groups are not well-defined. Indeed, in \cite{HZ} and \cite{HZ2} Zhang and myself proved that $h^{0,1}_{\bar\partial}$ can take different values on the same almost complex 4-manifold. In \cite{TT} Tardini and Tomassini proved the same was true for $h^{1,1}_{\bar \partial}$, although they were able to show that on compact almost Hermitian 4-manifolds $h^{1,1}_{\bar\partial}$ has a lower bound of $b_-$. In this paper we will give an upper bound of $b_- + 1$  and thereby obtain the result
\begin{repthm}{thm A}
   If $(X,J,\omega)$ is a compact almost Hermitian 4-manifold we have either $h^{1,1}_{\bar\partial} = b_-$ or $b_- +1$.
\end{repthm}

It is a widely-known result following from the classification of complex surfaces by Kodaira-Spencer (see \cite{BPV}) and later proven analytically by Buchdahl \cite{Buc} and Lamari \cite{Lam} that when $J$ is integrable on a compact 4-manifold, $b_1$ is even if and only if $(X,J)$ is K\"ahler, \textit{i.e.} if and only if $(X,J)$ admits an Hermitian metric $\omega$ such that $d\omega = 0$. This is equivalent to saying $h^{1,1}_{\bar\partial} = b_- +1$ if and only if $(X,J)$ is K\"ahler (see \textit{e.g.} \cite{DLZ}). The K\"ahler criterion using $b_1$ does not appear to extend to non-integrable almost complex structures, however we still expect $h^{1,1}_{\bar\partial}$ to detect almost K\"ahlerness (see Zhang's survey \cite{Zha0}). It was shown in \cite{TT}, that when an almost Hermitian metric is \textit{globally conformal to an almost K\"ahler metric} we have $h^{1,1}_{\bar\partial} = b_-+1$ and when the metric is \textit{strictly locally conformal to an almost K\"ahler metric} we have $h^{1,1}_{\bar\partial} = b_-$. 
We therefore ask the question
\begin{repq}{Q1}
On a compact almost Hermitian 4-manifold, does the value of $h^{1,1}_{\bar\partial}$ give a full description of whether an almost Hermitian metric is \textit{conformally almost K\"ahler}? Specifically, in the case when the metric is not \textit{locally conformally almost K\"ahler} (and thus also not \textit{globally conformally almost K\"ahler}) do we always have $h^{1,1}_{\bar\partial} = b_-$?
\end{repq}
Although the answer to this question is not yet known, we prove a similar result for the space of $d$-harmonic $(1,1)$-forms, which we denote $h_{d}^{1,1}$.
\begin{repthm}{d harmonic}
   On a compact almost Hermitian 4-manifold $(X,J,\omega)$, $h^{1,1}_{d}=b_- +1$ if $\omega$ is in the conformal class of an almost K\"ahler metric, otherwise $h_{d}^{1,1} = b_-$.
\end{repthm}

At the end of Section \ref{dbar sec} we calculate the value of $h^{1,1}_{\bar\partial}$ for a large family of almost Hermitian structures on the Kodaira-Thurston manifold and show that, at least for this example, Question \ref{Q1} is answered positively.

In \cite{PT} Piovani and Tomassini prove that on a compact Hermitian 4-manifold $h^{1,1}_{BC}$ can only ever be $b_-$ or $b_- +1$. They ask whether both of these values can be attained by some choice of metric, a question we will answer in this paper with the following theorem.
\begin{repthm}{thm B}
Given any compact almost Hermitian 4-manifold $(X,J,\omega)$, we have $h^{1,1}_{BC}=b_- +1$.
\end{repthm}
We also perform a calculation of $h^{2,1}_{BC}$ and $h^{1,2}_{BC}$ for a family of almost Hermitian structures on the Kodaira-Thurston manifold using the method developed in \cite{HZ} for turning PDEs into a collection of ODE and number theory problems. From these results we conclude the following
\begin{repthm}{metric inv}
   On a compact almost Hermitian 4-manifold, when $(p,q) = (0,0)$, $(1,0)$, $(0,1)$, $(2,0)$, $(1,1)$, $(0,2)$ or $(2,2)$, $h^{p,q}_{BC}$ is metric independent, but for $(p,q) = (2,1)$ and $(1,2)$ there exist examples for which $h^{p,q}_{BC}$ does vary with the metric.  
\end{repthm}
In addition to this, by building on some results of Chen and Zhang in \cite{CZ} we will prove that $h^{p,0}_{BC} = h^{0,p}_{BC}$ are birational invariants for all values of $p$.
\begin{repthm}{thm C}
Let $u:X\rightarrow Y$ be a degree one pseudoholomorphic map between compact almost complex 4-manifolds. Then $h^{p,0}_{BC}(X) = h^{p,0}_{BC}(Y)$ for any $p\in \{0,1,2\}$.
\end{repthm}

In the papers \cite{HZ,HZ2} Weiyi Zhang and myself present a calculation of $h^{0,1}_{\bar\partial}$ on the Kodaira-Thurston manifold, achieved through the introduction of a method for decomposing smooth functions which proved useful for solving linear PDEs. The key idea was to view the manifold as a torus bundle over $S^1$, thereby allowing for a Fourier expansion on each fibre. Further information can then be gained by considering the behaviour of the Fourier coefficients when travelling around the base space. In the last two sections of this paper we will show how the techniques used to decompose functions on the Kodaira-Thurston manifold can be applied to any torus bundle $M$ over $S^1$ given by $\mathbb{R}^{n+1}$ with points identified by 
\begin{equation*}
 \begin{pmatrix} t \\ \mathbf x \end{pmatrix} \sim
 \begin{pmatrix} t \\ \mathbf x + \mathbf \eta \end{pmatrix} \quad \text{ and } \quad  \begin{pmatrix} t \\ \mathbf x \end{pmatrix} \sim
 \begin{pmatrix} t + \xi \\ A^{\xi}\mathbf x \end{pmatrix} 
\end{equation*}
for all $\xi \in \mathbb{Z}$, $\mathbf \eta \in \mathbb{Z}^n$, with $A \in GL_n(\mathbb{Z})$.
By partitioning $\mathbb{Z}^n$ into its orbits under the action of the group generated by $A$ and separating the finite orbits from the infinite orbits, we produce the following decomposition of the space of $L^2$ functions on $M$. Here $\orb_{\mathbf y}$ denotes the orbit containing the element $\mathbf y \in \mathbb{Z}^n$.
\begin{repthm}{L2 decomp}
   The space of $L^2$ functions on $M$ decomposes in the following way.
  $$L^2(M) = \left(\widehat{\bigoplus_{\substack{\orb_{\mathbf y}
  \in \mathcal{O}\\ \abs{\orb_{\mathbf y}}=\infty}}}\mathcal{H}_{\mathbf y}\right) \oplus  \left(\widehat{\bigoplus_{\substack{\orb_{\mathbf y}
  \in \mathcal{O}\\ \abs{\orb_{\mathbf y}}= N<\infty}}}\widehat{\bigoplus_{t_0 \in \mathbb{Z}}}\mathcal{H}_{t_0, \mathbf y}\right), $$
  where
  $$\mathcal{H}_{\mathbf y} = \left\{ \sum_{\xi\in\mathbb{Z}} f(t+\xi)e^{2\pi i \mathbf y \cdot A^\xi \mathbf x} \,\middle|\, f \in L^2(\mathbb{R})  \right\}  $$
  and
  $$\mathcal{H}_{t_0, \mathbf y} = \left\{ C e^{2\pi i \frac{t_0 t}{N}}\sum_{\xi = 0}^{N-1} e^{2\pi i\left(\frac{t_0 \xi}{N}+\mathbf y \cdot A^{\xi} \mathbf x\right)}\, \middle|\, C\in \mathbb{C} \,\right\}. $$
  Here $\hat \oplus$ denotes the direct sum followed by the closure with respect to the $L^2$ norm.
\end{repthm}
Projection onto each of the components of this decomposition is given by the maps $\mathcal{F}_{\mathbf y}: L^2(M) \rightarrow \mathcal{H}_{\mathbf y}$
$$\mathcal{F}_{\mathbf y}(f)(t)=\int_{[0,1]^n} f(t,\mathbf x) e^{-2\pi i \mathbf y \cdot \mathbf x} \, d\mathbf{x} $$
and $\mathcal{G}_{t_0, \mathbf y}: L^2(M) \rightarrow \mathcal{H}_{t_0, \mathbf y}$
$$\mathcal{G}_{t_0,\mathbf y}(f) = \frac 1N \int_0^N \mathcal{F}_{\mathbf y}(f)(t)e^{-\frac{2\pi i t_0 t}{N}}\, dt .$$

In \cite{HZ, HZ2} this decomposition always leads to solving a combination of ODEs and lattice counting problems. In the last section we demonstrate how this decomposition could also reduce the PDEs deriving from the calculation of $h^{0,1}_{\bar\partial}$ to a recurrence relation problem.

\textit{Acknowledgements.} The author would like to thank Riccardo Piovani, Nicoletta Tardini and Adriano Tomassini whose recent work inspired the first half of this paper, in addition to Weiyi Zhang whose advice has been invaluable. The author's research is funded through the EPSRC Doctoral Training Partnership award.

\section{Preliminary results}
In this section we will recall some important facts about almost Hermitian manifolds which will be useful for proving the results of this paper in later sections. Let $(X,J)$ be an almost complex manifold. The existence of the almost complex structure $J$ induces a decomposition of the space of complex valued $k$-forms $\mathcal{A}_{\mathbb{C}}
^k$ into spaces of $(p,q)$-forms
$$\mathcal{A}_{\mathbb{C}}^k = \bigoplus_{p+q = k} \mathcal{A}^{p,q}. $$
This in turn leads to a decomposition of the exterior derivative $d:\mathcal{A}^k \rightarrow \mathcal{A}^{k+1}$ into the sum of 4 components
$$d = \mu + \partial + \bar\partial + \bar\mu $$
which change the bidegree of a $(p,q)$-forms by $(+2,-1), (+1,0), (0,+1)$ and $(-1,+2)$ respectively. We say that the almost complex structure $J$ is integrable when $\mu = \bar\mu = 0$, in which case it arises from a complex structure on $X$.
Given an almost Hermitian metric we define the $d$, $\bar\partial$ and Bott-Chern laplacians by
$$\Delta_{d} = dd^* + d^* d $$
$$\Delta_{\bar\partial} = \bar \partial \bar{\partial}^* + \bar{\partial}^* \bar\partial $$
$$ \Delta_{BC} = \partial\bar\partial\bar{\partial}^*\partial^* + \bar{\partial}^* \partial^* \partial \bar\partial +  \partial^* \bar\partial \bar{\partial}^* \partial + \bar{\partial}^* \partial \partial^* \bar\partial + \partial^*\partial + \bar{\partial}^*\bar\partial $$
along with the spaces of harmonic forms $\mathcal{H}_{d}^k = \ker \Delta_d |_{\mathcal{A}^k}$, $\mathcal{H}^{p,q}_{\bar\partial} = \ker \Delta_{\bar\partial}|_{\mathcal{A}^{p,q}}$ and  $\mathcal{H}^{p,q}_{BC} = \ker \Delta_{BC}|_{\mathcal{A}^{p,q}}$.
Here we define the adjoints of $\partial$ and $\bar\partial$ to be $\partial^* = -* \bar\partial *$ and $\bar{\partial}^* = - * \partial *$ where $*$ denotes the Hodge star operator. The dimensions of the spaces of $\bar\partial$ and Bott-Chern harmonic $(p,q)$-forms are denoted by $h_{\bar\partial}^{p,q}$ and $h_{BC}^{p,q}$. 

On a compact manifold the property of a general differential form $s$ being $d$, $\bar\partial$ or Bott-Chern harmonic can be equated to a collection of conditions as follows
$$\Delta_{d}s = 0 \iff \begin{cases}d s = 0\\ d * s =0 
\end{cases} $$ 
\begin{equation}\label{Delta equiv}
    \Delta_{\bar\partial}s = 0 \iff \begin{cases}\bar\partial s = 0\\ \partial * s =0 
\end{cases} \quad\quad\quad \Delta_{BC}s = 0 \iff \begin{cases} \partial s = 0 \\ \bar\partial s = 0\\ \partial \bar\partial * s =0 \end{cases}.
\end{equation} 

At this point it should be noted that the existence of an almost Hermitian metric $g$ is equivalent to the existence of a compatible $(1,1)$-form $\omega$ called the fundamental form, one being derived from the other by the formula $g(\cdot, \cdot) = \omega(\cdot, J \cdot)$. Consequently in this paper, as in many others, we will often refer to $\omega$ as if it were the corresponding almost Hermitian metric.   
If on a compact almost Hermitian 4-manifold, $\omega$ is Gauduchon, \textit{i.e.} $\partial \bar\partial \omega = 0$, then a result of Tardini and Tomassini \cite{TT} tells us that $\mathcal{H}^{1,1}_{\bar\partial}$ can be characterised by
\begin{equation}\label{charac 1}
    \mathcal{H}^{1,1}_{\bar\partial} = \{a\omega + \gamma \, | \, a \in \mathbb{C}, *\gamma = - \gamma, id^c \gamma = a d \omega\}
\end{equation}
whilst a result of Piovani and Tomassini \cite{PT} tells us that $\mathcal{H}^{1,1}_{BC}$ can be characterised by
\begin{equation}\label{charac 2}
\mathcal{H}^{1,1}_{BC} = \{ a\omega - \gamma \, | \, a \in \mathbb{C},  *\gamma = -\gamma, d\gamma = a d\omega \}.
\end{equation}
Here we define $d^c := J^{-1} d J$ with $J$ acting on a $(p,q)$-form as multiplication by $i^{p-q}$.

For any two conformal metrics $\omega = f \tilde\omega$ on a 4-manifold, the two resulting Hodge stars differ by $*_{\omega} = f^{2-p-q}*_{\tilde\omega}$ when acting on a $(p,q)$-form. From \eqref{Delta equiv} we can see that this means $\mathcal{H}^{1,1}_{\bar\partial}$ and $\mathcal{H}_{BC}^{1,1}$ are conformally invariant. Therefore, since a result of Gauduchon \cite{Gau} states that every conformal class contains a Gauduchon metric, we can apply \eqref{charac 1} and \eqref{charac 2} given any almost Hermitian metric by finding the Gauduchon metric to which it is conformal.   

For any almost Hermitian metric $\omega$ we have the property that
$$d \omega = \alpha \wedge \omega $$
for some 1-form $\alpha$. This comes as a consequence of the well known fact that the map $L^k:\mathcal{A}^{n-k} \rightarrow \mathcal{A}^{n+k}$ given by $s \mapsto s \wedge\omega^k  $ is bijective. Furthermore, $\alpha$ is an exact form if and only if $\omega$ is \textit{globally conformal to an almost K\"ahler metric}, while $\alpha$ is closed if and only if $\omega$ is \textit{locally conformal to an almost K\"ahler metric}.

\section{$\bar{\partial}$-harmonic $(1,1)$-forms}\label{dbar sec}
From the characterisation \eqref{charac 1} of $\mathcal{H}^{1,1}_{\bar\partial}$ we conclude the following.
\begin{thm}\label{thm A}
  If $(X,J,\omega)$ is a compact almost Hermitian 4-manifold we have either $h_{\bar\partial}^{1,1} = b_-$ or $b_- + 1$.
  \begin{proof}
    Since $h_{\bar\partial}^{1,1}$ is a conformal invariant we assume without loss of generality that $\omega$ is a Gauduchon metric. 
    
    From \eqref{Delta equiv} we obtain the inclusion
    $$\mathcal{H}_g^- \subseteq \mathcal{H}_{\bar\partial}^{1,1}  $$
    where $\mathcal{H}_g^-$ denotes the space of $d$-harmonic anti-self-dual $(1,1)$-forms. When this inclusion is an equality then clearly we have $h_{\bar\partial}^{1,1} = b_-$.
    Suppose instead that $\mathcal{H}_{\bar\partial}^{1,1}$ has some element $a_0 \omega + \gamma_0$ which is not in $\mathcal{H}_g^-$. Here $a_0$ is a constant and $\gamma_0$ is an anti-self-dual form satisfying $id^c \gamma_0 = a_0 d\omega$. Note that $a_0$ cannot be zero, as that would leave us with a $d$-harmonic anti-self-dual form.   
     A general element of $\mathcal{H}_{\bar\partial}^{1,1}$ given by $a \omega + \gamma$ can then be rewritten as an element of $\mathcal{H}_g^-$ plus a multiple of the single additional element $a_0 \omega + \gamma_0 $
    $$a\omega + \gamma =  \frac {a}{a_0} \left(a_0\omega + \gamma_0\right) + \frac{1}{a_0}\left(a_0 \gamma - a \gamma_0\right),$$
    thus giving us $h_{\bar\partial}^{1,1} = b_-+1$.
    
    To see that the anti-self-dual form $a_0\gamma - a \gamma_0$ is indeed $d$-harmonic, first note that $d^c(a_0\gamma - a \gamma_0) = a_0 d^c \gamma - a d^c \gamma_0 =0$. Then, since $d^c = J^{-1} d J $ and $J$ is the identity when acting on $(1,1)$-forms, it follows that $d(a_0 \gamma - a \gamma_0)=0$. As our form is anti-self-dual we therefore also have $d * (a_0 \gamma - a \gamma_0) = 0$.   
  \end{proof}
\end{thm}
\begin{cor}
    If $(X,J,\omega)$ is a compact almost Hermitian 4-manifold where we assume $\omega$ is Gauduchon, then $h^{1,1}_{\bar\partial}=b_- + 1$ if and only if there exists an anti-self-dual $(1,1)$-form $\gamma$ satisfying the equation
    \begin{equation}\label{delbar eq}
        id^c \gamma= d \omega.
    \end{equation}
    \begin{proof}
          If such a $\gamma$ exists then $ \omega + \gamma$ is $\bar\partial$-harmonic, along with $b_-$ many linearly independent elements of $\mathcal{H}_{g}^-$, therefore $h^{1,1}_{\bar\partial}=b_- +1$.
          
          Conversely, if $h^{1,1}_{\bar\partial}=b_- +1$, then there must be some form in $\mathcal{H}_{\bar\partial}^{1,1}$ other than those contained in $\mathcal{H}_{g}^{-}$, \textit{i.e.} a form which can be written as $a_0 \omega + \gamma_0$ with $a_0 \neq 0$ such that $id^c \gamma_0 = a_0 d\omega$. Thus $\gamma =\frac{1}{a_0}\gamma_0$ gives us the desired solution. 
    \end{proof}
\end{cor}

In \cite{DLZ}, Draghici, Li and Zhang prove that, for integrable almost complex manifolds $(X,J)$, $h_{\bar\partial}^{1,1}$ takes the value $b_- + 1$ when $(X,J)$ is K\"ahler and otherwise takes the value $b_-$. Partially extending this result to non-integrable manifolds, in \cite{TT} it was proven that if a compact almost Hermitian 4-manifold $(X,J,\omega)$ is \textit{globally conformally almost K\"ahler} then $h_{\bar\partial}^{1,1} = b_- +1$, whereas if $(X,J,\omega)$ is \textit{strictly locally conformally almost K\"ahler} then $h_{\bar\partial}^{1,1} = b_-$. We therefore ask the question
\begin{q}\label{Q1}
    On a compact almost Hermitian 4-manifold, does the value of $h_{\bar\partial}^{1,1}$ give a full description of whether an almost Hermitian metric is \textit{conformally almost K\"ahler}? Specifically, in the case when the metric is not \textit{locally conformally almost K\"ahler} (and thus also not \textit{globally conformally almost K\"ahler}) do we have $h^{1,1}_{\bar\partial} = b_-$?
\end{q}

Although the answer to this is not known, we can prove a similar result for the dimension of the space of $d$-harmonic $(1,1)$-forms, which we will denote by $h_d^{1,1}$.

\begin{thm}\label{d harmonic}
  On a compact almost Hermitian 4-manifold $(X,J,\omega)$, $h_d^{1,1}=b_- + 1$ if $\omega$ is in the conformal class of an almost K\"ahler metric, otherwise $h_d^{1,1}=b_-$.
  \begin{proof}
    As in the proof of the previous theorem, we use the fact that $h_d^{1,1}$ is a conformal invariant and thereby assume $\omega$ is a Gauduchon metric. Furthermore, all almost K\"ahler metrics are Gauduchon, so the conformal class of $\omega$ contains an almost K\"ahler metric if and only if $\omega$ is almost K\"ahler itself.   
    
    On compact manifolds we know a differential form $s$ is $d$-harmonic if and only if
    $$ds = 0 \quad \quad \quad d * s = 0. $$
    From this we can see that the Hodge star maps $d$-harmonic forms to $d$-harmonic forms, meaning that if some $(1,1)$-form $s$ is in $\mathcal{H}_d^{1,1}$ so too are its self-dual and anti-self-dual components, $\frac 12 (s + *s)$ and $\frac 12 (s - *s)$. Furthermore, we have the inclusion 
    $$\mathcal{H}_d^{1,1} \subseteq \mathcal{H}_{\bar\partial}^{1,1} $$
    and so from \eqref{charac 1} we know we can write any $d$-harmonic $(1,1)$-form as $a \omega + \gamma$ with $a \in \mathbb{C}$ a constant and $\gamma$ an anti-self-dual form. But the self-dual component of this is only harmonic if $d\omega = 0$ or $a = 0$ and so either $\omega$ is almost K\"ahler and we have $h_{d}^{1,1} = b_- +1$ or all $d$-harmonic $(1,1)$-forms are anti-self-dual and we have $h_{d}^{1,1} = b_-$.
  \end{proof}
\end{thm}
From this result we see that the above question is equivalent to asking whether $h^{1,1}_d$ and $h_{\bar\partial}^{1,1}$ are always equal on compact Hermitian 4-manifolds. 

We conclude this section with a calculation of $h^{1,1}_{\bar \partial}$ for a large family of almost complex structures and compatible metrics. In doing so we will see that, at least for this family of almost Hermitian structures, Question \ref{Q1} has a positive answer.

\begin{ex}\label{ex 1}
Briefly we recall the definition of the Kodaira-Thurston manifold $KT^4 = \Gamma \backslash G$ as the group $G = \mathbb{R}\times Nil^3$ modulo the subgroup $\Gamma$ of elements with integer valued entries, acting on $G$ by left multiplication. This is equivalent to $\mathbb{R}^4 /\sim$ with points identified by the equivalence relation
$$\begin{pmatrix} t \\ x \\ y \\ z \end{pmatrix} \sim
\begin{pmatrix} t + t_0 \\ x + x_0 \\ y + y_0 \\ z + z_0 + t_0 y \end{pmatrix}$$
for all $t_0, x_0, y_0, z_0 \in \mathbb{Z}$. 

$\mathbb{R}\times \Nil^3$ has a smooth global frame given by
$$\der t \quad \quad \der x \quad \quad \der y + t \der z \quad \quad \der z$$
which descends to a global frame for $KT^4$ since all of the above vector fields are invariant under the action of $\Gamma$.
We can define an almost complex structure acting on this frame by the matrix
$$J_{a,b} = \begin{pmatrix}
        0 & -1 & 0 &  0\\
        1 &  0 & 0 &  0\\
        0 &  0 & a &  b\\
        0 &  0 & c &  -a
        \end{pmatrix},$$
with $a, b\in \mathbb R$, $b \neq 0$ and $c = -\frac{a^{2}+1}{b}$.
A pair of vector fields
$$V_{1} = \frac 12 \left(\der{x}-i\der{t}\right)\quad \mathrm{\&}\quad V_{2} = \frac 12 \left(\left(\der{y}+t\der{z}\right)-\frac{a-i}{b}\der{z}\right) $$
can then be defined, spanning $T^{1,0}_{p}KT^4$ at every point $p\in KT^4$. Their dual $(1,0)$-forms are given by
$$\phi_{1} = \dif{x}+i\dif{t}\quad\mathrm{\&}\quad \phi_{2} = (1-ai)\dif{y}-ib(\dif{z}-t\dif{y}). $$
These forms satisfy the structure equations
$$d\phi^1 = 0 \quad \quad d\phi^2 = \frac b4 \left( \phi^{12}+\phi^{1\bar2} + \phi^{2\bar1} - \phi^{\bar1 \bar2}\right) $$
with $\phi^{i\bar j}$ used here as shorthand for $\phi^i \wedge \bar{\phi^j}$. From this we can see that $J_{a,b}$ is a non-integrable almost complex structure, namely we have $\bar \mu \phi^2 = - \frac b4 \phi^{\bar 1 \bar 2} \neq 0$.

Now it only remains for us to choose a family of almost Hermitian metrics 
$$\omega_{w} = i \left((1 + \abs{w}^2)\phi^{1\bar 1} - w \phi^{1\bar2} - \bar w \phi^{2\bar1} + \phi^{2\bar2}\right)$$ 
varying over some complex number $w \in \mathbb{C}$, defined such that $V_1 + w V_2$ and $V_2$ form a unitary basis on $T^{1,0}_p KT^4$. 
Using the structure equations we can calculate 
$$d \omega_{w} = iw\frac{b}{2} \phi^1 \wedge \bar{\phi^ 1} \wedge (\phi^2 + \bar{\phi^2}) $$
from which we see firstly that $\omega_w$ is an almost K\"ahler metric if and only if $w = 0$ and secondly that 
$$\partial \bar \partial \omega_w = iw\frac{b}{2} \phi^1 \wedge \bar{\phi^1} \wedge \partial\bar{\phi^2}  = 0$$
and thus $\omega_w$ is Gauduchon for all $w$.
Furthermore, we can write 
$$d\omega_w = \alpha_w \wedge \omega_w $$
with 
$$\alpha_w =   b \left( -w^2 \phi^1 - \bar{w}^2 \phi^{\bar 1} + w \phi^2 + \bar w \phi^{\bar 2} \right),$$
$$d \alpha_w = \frac{b^2}{4}(w - \bar w)(\phi^{12}+\phi^{1\bar2} + \phi^{2\bar1} - \phi^{\bar1 \bar2}).$$
$\omega_w$ is therefore \textit{globally conformally almost K\"ahler} only when $w = 0$ and \textit{locally conformally K\"ahler} only when $w$ takes real values.

Finding $h_{\bar\partial}^{1,1}$ then amounts to asking whether there exists an anti-self-dual $\gamma$ solving 
$$id^c \gamma = d \omega_w. $$
Since $J$ is the identity on $(1,1)$-forms this is equivalent to 
$$iJ^{-1}d\gamma =  d\omega_w. $$
If such a $\gamma$ exists that would mean
\begin{align*}
    J d\omega_w &= w\frac{b}{2} \phi^1 \wedge \phi^{\bar 1} \wedge (\phi^2 - \phi^{\bar 2})  \\
    &= -2wb \, dx \wedge dt \wedge ((a-t)dy + b dz)
\end{align*}
is an exact 3-form, but consider the closed submanifold given by $y=0$. The pullback of $J \omega_w$ onto this submanifold is $-2wb^2 \, dx\wedge dt \wedge dz $, which by Stokes' theorem cannot be exact since its integral over the submanifold is non-zero, the only exception to this being when $w = 0$. Thus, in all the cases when $\omega_w$ is not globally almost K\"ahler, there is no solution to \eqref{delbar eq} and so $h^{1,1}_{\bar\partial} = b_-$ = 2. 
\end{ex}

\section{Bott-Chern harmonic forms}
In this section we will give a collection of results which together will give a full description of when $h^{p,q}_{BC}$ is or is not metric independent for compact 4-manifolds.

For many values of $(p,q)$ proving the metric invariance of $h^{p,q}_{BC}$ is a relatively trivial affair and so we will not spend too long on these cases.
\begin{lem}
  On any compact almost Hermitian 4-manifold $h^{p,q}_{BC}$ is metric independent when $(p,q)$ is equal to $(2,0), (0,2), (1,0), (0,1), (0,0)$ or $(2,2)$. 
  \begin{proof}
        Bott-Chern harmonic $(0,0)$-forms are always just the constant functions, since $\Delta_{BC}$ is elliptic. Similarly Bott-Chern harmonic $(2,2)$-forms are just constant functions times the volume form so although $\mathcal{H}^{2,2}_{BC}$ might change with the metric, $h^{2,2}_{BC}$ does not.
        
        For the remaining cases recall that a $(p,q)$-form $s$ is Bott-Chern harmonic if and only if it satisfies the three conditions
        $$\partial s = 0 \quad \quad \bar\partial s = 0 \quad \quad \partial \bar\partial * s = 0.$$
        When $(p,q) = (2,0), (0,2), (1,0)$ or $(0,1)$ the third condition is always true leaving behind the first two conditions which do not depend on the metric.
  \end{proof}
\end{lem}

The more interesting cases are those when $(p,q) = (1,1), (2,1)$ and $(1,2)$. We start with the case of $\mathcal{H}^{1,1}_{BC}$. From the characterisation \eqref{charac 2} of $\mathcal{H}^{1,1}_{BC}$, in \cite{PT} it is deduced that $h_{BC}^{1,1}$ is either $b_- +1$ or $b_-$, with the two cases corresponding, respectively, to the existence or non-existence of an anti-self-dual solution $\gamma$ to the equation
 \begin{equation}\label{BC cond}
     d\gamma = d\omega.
 \end{equation}
Here $\omega$ is a Gauduchon metric conformal to the chosen Hermitian metric.

It turns out that solutions to the above equation can be found by making use of the Hodge decomposition
$$\mathcal{A}^k =  d\mathcal{A}^{k-1} \oplus \mathcal{H}_{d}^{k} \oplus d^*\mathcal{A}^{k+1}.$$
\begin{thm}\label{thm B}
Given any compact almost Hermitian 4-manifold $(X,J,\omega)$, we have $h^{1,1}_{BC} = b_- +1$.
\begin{proof}
From the conformal invariance of $h_{BC}^{1,1}$ we may assume without losing generality that $\omega$ is Gauduchon. Then taking the Hodge decomposition we can write
$$\omega = d\alpha + h + d^* \beta$$
for some $\alpha \in \mathcal{A}^1, h\in\mathcal{H}_{d}^{2}$ and $\beta \in \mathcal{A}^3$. By defining a 2-form
$$\gamma = d *\beta + d^* \beta  $$
we have
$$d\omega = dd^* \beta = d\gamma$$
and thus $\gamma$ is a solution to \eqref{BC cond}. 

It only remains to show that $\gamma$ is anti-self-dual. 
Using the definition of $d^*$ along with the fact that the square of the Hodge star when applied to a $k$-form is given by $*^2 = (-1)^k$, we can see that
\begin{align*}
    *\gamma &= *d*\beta - **d*\beta\\
    &= -d^* \beta - d * \beta\\
    &= -\gamma.
\end{align*}
We therefore find that 
$$\mathcal{H}_{BC}^{1,1} = \mathcal{H}^-_g  \oplus  \mathbb{C}\langle \omega - \gamma \rangle $$
and so $h_{BC}^{1,1}$ is always $b_- + 1$.
\end{proof}
\end{thm}

We will now use the following example to show that $h^{2,1}_{BC}$ and $h^{1,2}_{BC}$ may, in general, depend on the choice of almost Hermitian metric.

\subsection{Calculating $h_{BC}^{2,1}$ and $h_{BC}^{1,2}$ on the Kodaira-Thurston manifold}\label{example KT4}
For this example we again consider the Kodaira-Thurston manifold, with the same almost complex structure as in Example \ref{ex 1}. The metric we will be using is given by
$$\omega_{\rho} = i \left( \phi^{1  \bar1} + \rho \phi^{2 \bar 2}  \right) $$
such that $V^1, \bar V^1 , \frac{1}{\sqrt \rho}V^2 $ and $\frac{1}{\sqrt \rho}\bar V^2$ form an orthonormal basis.
This is essentially the same metric as was used in \cite{HZ2} and in fact what follows is a more general, completed version of a calculation in \cite{PT}.

\begin{ex} 
Let a general $(2,1)$-form be given by $f \phi^{12\bar1} + g \phi^{12\bar2}$. Then from the conditions $\bar\partial s = 0$ and $\partial \bar\partial * s = 0$ we see that $s \in \mathcal{H}^{2,1}_{BC}$ if and only if the following PDEs hold.
\begin{equation}\label{2,1}
    \begin{cases}
\rho V_1 \bar V_1 (f) +  V_2 \bar V_1 (g) - \frac b4 \rho V_1 (f) + \frac b4 \rho \bar V_1 (f) - \frac b4 \bar V_2 (g) - \frac{b^2}{8} \rho f= 0\\
\rho V_1 \bar V_2 (f) + V_2 \bar V_2 (g) + \frac b4 \rho V_2 (f) = 0 \\
\bar V_1 (g) - \bar V_2 (f) = 0
\end{cases}
\end{equation}
Using the same method as in $\cite{HZ2}$ we can perform a Fourier expansion with respect to $x,y$ and $z$ to simplify the above equations. We will write 
$$f(t,x,y,z) = \sum_{l,m,n}\mathcal{F}_{l,m,n}(f)(t)e^{2\pi i (lx + my + nz)} $$
where 
$$\mathcal{F}_{l,m,n}(f)(t)=\int_{[0,1]^3} f(t,x,y,z) e^{-2\pi i (lx+my+nz)} \, dx \, dy\, dz.$$
Applying a Fourier expansion to the second and third PDEs we obtain the ODE system 
\begin{equation}\label{ODE 2,1}
    \frac{d}{dt} \begin{pmatrix}
\mathcal{F}_{l,m,n}(f) \\ \mathcal{F}_{l,m,n}(g)
\end{pmatrix} = 2\pi \left[ \begin{pmatrix}
0 & \frac n\rho \\ n & 0
\end{pmatrix}t +  \begin{pmatrix}
l - \frac{b}{4\pi} i & \frac 1\rho \left(m - n\frac{a-i}{b}\right) \\
m-n\frac{a+i}{b} & -l
\end{pmatrix} \right]\begin{pmatrix}
\mathcal{F}_{l,m,n}(f) \\ \mathcal{F}_{l,m,n}(g)
\end{pmatrix} 
\end{equation}
for every $l,m,n \in \mathbb{Z}$.
The ODE given by expanding our first PDE can be derived from the above ODE system and so adds no new information.

As was proven in \cite{HZ2}, the solutions to $\eqref{2,1}$ can be split into two cases:

Firstly if two smooth functions $\mathcal{F}_{l,m,n}(f)$, $\mathcal{F}_{l,m,n}(g) \in C^\infty (\mathbb{R})$ satisfy the ODE $\eqref{ODE 2,1}$ with $n \neq 0$ and $0 \leq m < \abs{n}$ then we have a solution to \eqref{2,1} given by
$$f = \sum_{\xi \in \mathbb{Z}}\mathcal{F}_{l,m,n}(f)(t+\xi)e^{2\pi i (lx+(m+n\xi)y +nz)} $$
$$g = \sum_{\xi \in \mathbb{Z}}\mathcal{F}_{l,m,n}(g)(t+\xi)e^{2\pi i (lx+(m+n\xi)y +nz)} $$
if and only if the two functions $\mathcal{F}_{l,m,n}(f)$ and $\mathcal{F}_{l,m,n}(g)$ are Schwartz.

Secondly, if $\mathcal{F}_{l,m,n}(f)$, $\mathcal{F}_{l,m,n}(g) \in C^\infty (\mathbb{R})$ satisfy the ODE $\eqref{ODE 2,1}$ with $n = 0$ then we have a solution to \eqref{2,1} given by
$$f = \mathcal{F}_{l,m,0}(f)e^{2\pi i (lx + my)} $$
$$g = \mathcal{F}_{l,m,0}(g)e^{2\pi i (lx+my)} $$
if and only if the two functions $\mathcal{F}_{l,m,n}(f)$ and $\mathcal{F}_{l,m,n}(g)$ are periodic with a period of 1.

Finding solutions in the first case amounts to solving a Stokes phenomenon problem. This can be tricky to do in general, but this problem has been solved for the ODE \eqref{ODE 2,1} in Theorem 3.1 of \cite{HZ}. It turns out we have a solution for all $0\leq m < \abs{n}$ whenever $l = 0$ and $n$ satisfies
$$64 \pi^2 n^2 - 64  \pi n u b^2 \sqrt{\rho} - b^4 \rho =0 $$
for some negative integer $u$. Or equivalently, if we set $d = \frac{b}{8\pi}$,
$$n^2 - 64\pi nud^2 \sqrt{\rho}- 64 \pi^2 d^4 \rho= 0. $$
Note that if $d$ and $\rho$ are both rational this case gives us no solutions as $\pi$ is transcendental.

For the second case, since we are working with periodic functions, we can take another Fourier expansion with respect to $t$, writing
$$\mathcal{G}_{k,l,m,0}(f) = \int_0^1 \mathcal{F}_{l,m,0}(f)(t) e^{-2\pi i kt } dt. $$
Applying this expansion to $\eqref{ODE 2,1}$ we obtain the equations
$$\rho\left(l - ik - \frac{b}{4\pi} i \right) \mathcal{G}_{k,l,m,0}(f) + m \mathcal{G}_{k,l,m,0}(g) = 0 $$
$$m \mathcal{G}_{k,l,m,0}(f) = (l+ik) \mathcal{G}_{k,l,m,0}(g). $$
This can be solved directly to find the solution
$$s =  \phi^{12\bar 2}$$
when $k = 0$, and the solution
$$s =  i k e^{2\pi i(kt+my)}\phi^{12\bar 1}+ m e^{2\pi i(kt+my)}\phi^{12\bar 2} $$
when $k\neq 0$ and $k,m\in \mathbb{Z}$ satisfy
$$\frac{m^2}{\rho} + (k+d)^2 = d^2. $$
Here we again set $d = \frac{b}{8\pi}$. Notice that when $d = 1$ and $\rho = 1$ we have 4 solutions given by $(k,m) = (-1,1), (-1,-1) (-2,0)$ and $(0,0)$, however when when we take $\rho = \frac 12$, leaving $d$ unchanged, we only have the two solutions $(k,m) = (-1,0)$ and $(0,0)$. Therefore we conclude that on the Kodaira-Thurston manifold the value of $h^{2,1}_{BC}$ may depend on the choice of almost Hermitian metric. 
\end{ex}

\begin{ex}
Now let a general $(1,2)$-form be given by $f \phi^{1\bar1 \bar 2} + g \phi^{2\bar1 \bar2}$. Then from the conditions $\partial s = 0$ and $\partial \bar\partial * s = 0$ we see that $s \in \mathcal{H}^{1,2}_{BC}$ if and only if the following PDEs hold.
\begin{equation}\label{1,2}
    \begin{cases}
\rho V_1 \bar V_1 (f) +  V_1 \bar V_2 (g) + \frac b4 \rho V_1 (f) - \frac b4 \rho \bar V_1 (f) -\frac b4 \bar V_2 (g) - \frac{b^2}{16} \rho f = 0\\
\rho V_2 \bar V_1 (f) + V_2 \bar V_2 (g) + \frac b4 \rho V_2 (f) = 0 \\
V_1 (g) -  V_2 (f) = 0
\end{cases}
\end{equation}
Applying the same Fourier expansion as before, the second and third equations give us the ODE system
\begin{equation*}
    \frac{d}{dt} \begin{pmatrix}
\mathcal{F}_{l,m,n}(f) \\ \mathcal{F}_{l,m,n}(g)
\end{pmatrix} = 2\pi \left[ \begin{pmatrix}
0 & \frac{n}{\rho} \\ n & 0
\end{pmatrix}t +  \begin{pmatrix}
-l + \frac{b}{4\pi} i & -\frac{1}{\rho}\left(m - n\frac{a+i}{b}\right) \\
-m + n\frac{a-i}{b} & l 
\end{pmatrix} \right]\begin{pmatrix}
\mathcal{F}_{l,m,n}(f) \\ \mathcal{F}_{l,m,n}(g)
\end{pmatrix}.
\end{equation*}

Again splitting the solutions into two cases we find that firstly we have a solution for all $n \neq 0$ and $0 \leq m< \abs{n}$ which satisfy 
$$n^2 - 64\pi nud^2 \sqrt{\rho}- 64 \pi^2 d^4 \rho= 0 $$
for some negative integer $u$.
Secondly, for $n = 0$ we have solutions 
$$s =  \phi^{2\bar 1 \bar 2}$$
and
$$s =  i k e^{2\pi i(kt+my)}\phi^{1\bar 1 \bar 2}- m e^{2\pi i(kt+my)}\phi^{2\bar 1 \bar 2} $$
for all $k,m\in \mathbb{Z}$, with $k \neq 0$, satisfying
$$\frac{m^2}{\rho} + (k-d)^2 = d^2. $$
From the above we see that for this family of almost Hermitian structures we have $h^{1,2}_{BC}= h^{2,1}_{BC}$ (although this need not always be the case). Thus the value of $h^{1,2}_{BC}$ may also depend on the choice of almost Hermitian metric.

Furthermore, when $\rho = 1$, the calculation of Theorem 4.1 in \cite{HZ} tells us that $h^{2,1}_{BC}$ and $h_{BC}^{1,2}$ here are both equal to $h_{\bar\partial}^{0,1}$ defined using the same family of almost complex structures $J_{a,b}$. In particular, $h_{BC}^{2,1}$ and $h_{BC}^{1,2}$ can both be made arbitrarily large by varying the value of $b$.   
\end{ex}

We can now bring the results of this section together into the following theorem.

\begin{thm}\label{metric inv}
  On a compact almost Hermitian 4-manifold, when $(p,q) = (0,0)$, $(1,0)$, $(0,1)$, $(2,0)$, $(1,1)$, $(0,2)$ or $(2,2)$, $h^{p,q}_{BC}$ is metric independent, but for $(p,q) = (2,1)$ and $(1,2)$ there exist examples for which $h^{p,q}_{BC}$ does vary with the metric.  
\end{thm}

\subsection{Birational invariance of $h^{p,0}_{BC}$}
It is known from Theorem 5.5 in \cite{CZ} that $h^{p,0}_{\bar{\partial}}$ is birationally invariant on compact 4-manifolds for any $p \in \{0,1,2\}$. This means that if we have a sequence of almost complex 4-manifolds $X_0, X_1, X_2 \dots, X_{k+1}$ along with a sequence of degree one pseudoholomorphic maps $u_0, \dots u_k$ such that $u_{2i-1}:X_{2i-1}\rightarrow X_{2i}$ and $u_{2i}:X_{2i+1}\rightarrow X_{2i}$ then $h^{p,0}_{\bar\partial}(X) = h^{p,0}_{\bar\partial}(Y)$. It turns out this result can be extended to show that the numbers $h^{p,0}_{BC}$ are also birational invariants.
\begin{thm}\label{thm C}
  Let $u:X\rightarrow Y$ be a degree one pseudoholomorphic map between compact almost complex 4-manifolds. Then $h^{p,0}_{BC}(X) = h^{p,0}_{BC}(Y)$ for any $p\in \{0,1,2\}$. 
  \begin{proof}
    From \cite{CZ} we know that the pullback with respect to $u$ describes a bijection $$u^*:\mathcal{H}^{p,0}_{\bar\partial}(Y)\rightarrow \mathcal{H}^{p,0}_{\bar\partial}(X).$$
    Restricting this to the forms $s\in \mathcal{H}^{p,0}_{\bar\partial}(Y)$ which satisfy $\partial s = 0$ gives us
    $$u^*:\mathcal{H}^{p,0}_{BC}(Y)\rightarrow \mathcal{H}^{p,0}_{BC}(X).$$
    The injectivity of this map follows directly from the injectivity of $u^*$ acting on $\mathcal{H}^{p,0}_{\bar\partial}(Y)$, so it only remains to prove surjectivity.
    
    Since $u^*$ is invertible when acting on $\mathcal{H}^{p,0}_{\bar\partial}(Y)$ we know that for any $s\in \mathcal{H}_{BC}^{p,0}(X)$ there is some $t\in \mathcal{H}^{p,0}_{\bar\partial}(Y)$ such that $u^* t = s$. By Theorem 1.5 in \cite{Zha} we know there is a finite set $Y_1\subset Y$ such that the restriction 
    $$u: X\backslash u^{-1}(Y_1)\rightarrow Y\backslash Y_1$$ is a diffeomorphism. This means we have
    $$t{\big |}_{X\backslash u^{-1}(Y_1)} = (u^{-1})^*s{\big |}_{Y\backslash Y_1} $$
    and so $\partial t = 0$ on $Y\backslash Y_1$. But since $t$ is smooth and $\overline{Y\backslash Y_1} = Y$, we must have $\partial t = 0$ on all of $Y$, thus $t\in \mathcal{H}^{p,0}_{BC}(Y)$ and $u^*{\big |}_{\mathcal{H}^{p,0}_{BC}(Y)}$ is surjective. 
  \end{proof}
\end{thm}
\begin{cor}
  $h^{0,p}_{BC}$ is a birational invariant on compact almost complex 4-manifolds for any $p = 0,1$ or $2$.
  \begin{proof}
    Recall that $s\in \mathcal{H}^{p,q}_{BC}$ if and only if the following conditions hold 
    $$\bar\partial s = 0 \quad \quad \partial s = 0 \quad \quad \partial \bar\partial * s = 0.$$
    If $s$ is either a $(p,0)$-form or a $(0,p)$-form for any $p = 0, 1$ or $2$ then the third condition is always true for reasons of bidegree. The remaining two conditions, when taken together, are unchanged by a conjugation of $s$. The corollary therefore follows simply from the fact that $\mathcal{H}_{BC}^{0,p} = \overline{\mathcal{H}_{BC}^{p,0}}$.
  \end{proof}
\end{cor}

\section{Harmonic Analysis on Torus bundles over $S^1$}

In this section we introduce a technique which may be used to simplify or solve certain linear PDEs on torus bundles over $S^1$. Special cases of this technique have already proven useful in the calculation of $h^{p,q}_{\bar\partial}$ on the Kodaira-Thurston manifold \cite{HZ}. We will start by first describing a decomposition of smooth functions. Then, by considering a specific example of calculating $h^{0,1}_{\bar\partial}$ on a torus bundle with Euclidean geometry, we will see how PDEs can be simplified through the application of this decomposition. In our example it will simplify to a recurrence relation.    

\subsection{Decomposition of functions}
Let $M$ be any $n$-torus bundle over $S^1$. This can be described as the mapping torus of an $n$-torus determined by a matrix $A\in GL_n(\mathbb{Z})$. In other words, $M$ is given by $\mathbb{R}^{n+1}$ with points identified by 

\begin{equation}
 \begin{pmatrix} t \\ \mathbf x \end{pmatrix} \sim
 \begin{pmatrix} t \\ \mathbf x + \mathbf \eta \end{pmatrix} \quad \text{ and } \quad  \begin{pmatrix} t \\ \mathbf x \end{pmatrix} \sim
 \begin{pmatrix} t + \xi \\ A^{\xi}\mathbf x \end{pmatrix} 
 \label{identification}
\end{equation}
for all $\xi \in \mathbb{Z}$, $\mathbf \eta \in \mathbb{Z}^n$. 

When $t$ is fixed, $\mathbf x$ describes a point on a torus. This means any smooth function $f\in C^\infty(M)$, when viewed as a function on $\mathbb{R}^{n+1}$ satisfying
  \begin{equation}
   f(t,\mathbf x ) = f(t,\mathbf x + \mathbf \eta ) \quad \text{ and } \quad f(t,\mathbf x) = f(t+\xi, A^{\xi}\mathbf{x}) 
   \label{periodicity}
  \end{equation}
can be decomposed into the Fourier series
$$f(t,\mathbf x) = \sum_{\mathbf{x}_0 \in \mathbb{Z}^n}\mathcal{F}_{\mathbf{x}_0}(f)(t)e^{2\pi i \mathbf{x}_0 \cdot \mathbf x} $$
where we define
$$\mathcal{F}_{\mathbf{x}_0}(f)(t)=\int_{[0,1]^n} f(t,\mathbf x) e^{-2\pi i \mathbf{x}_0 \cdot \mathbf x} \, d\mathbf{x} .$$
Here we have to be careful: notice that we have no guarantee that the summands $\mathcal{F}_{\mathbf{x}_0}(f)e^{2\pi i \mathbf{x}_0 \cdot \mathbf x} $ will satisfy the same condition \eqref{periodicity} as $f$, and so the summands are not themselves smooth functions on $M$. In particular, it is the second condition of \eqref{periodicity} that may fail. We do however have the following result.

\begin{prop}\label{function on M}
    A function $f\in C^\infty(\mathbb{R}^{n+1})$ satisfies \eqref{periodicity} if and only if it can be written as the Fourier series
    $$f(t, \mathbf{x})=\sum_{\mathbf{x}_0 \in \mathbb{Z}^n}\mathcal{F}_{\mathbf{x}_0}(f)(t)e^{2\pi i \mathbf{x}_0 \cdot \mathbf x} $$
    such that
    $$\mathcal{F}_{(A^T)^\xi \mathbf{x}_0}(t) = \mathcal{F}_{\mathbf{x}_0}(f)(t+\xi) $$
    for all $\xi \in\mathbb{Z}$.
    \begin{proof}
        It is clear that $f$ has a Fourier expansion if and only if it satisfies the first condition of \eqref{periodicity}. Taking the expansion of the second condition we see that
        $$\sum_{\mathbf{x}_0 \in \mathbb{Z}^n}\mathcal{F}_{\mathbf{x}_0}(f)(t)e^{2\pi i \mathbf{x}_0 \cdot \mathbf x}  = \sum_{\mathbf{x}_0 \in \mathbb{Z}^n}\mathcal{F}_{\mathbf{x}_0}(f)(t+\xi)e^{2\pi i \mathbf{x}_0 \cdot A^\xi \mathbf x}  $$
        or equivalently
        $$\sum_{\mathbf{x}_0 \in \mathbb{Z}^n}\mathcal{F}_{\mathbf{x}_0}(f)(t)e^{2\pi i \mathbf{x}_0 \cdot \mathbf x} =\sum_{\mathbf{x}_0 \in \mathbb{Z}^n}\mathcal{F}_{(A^T)^{-\xi}\mathbf{x}_0}(f)(t+\xi)e^{2\pi i \mathbf{x}_0 \cdot \mathbf x}  .$$
        By the uniqueness of Fourier coefficients, this is identical to requiring
        $$\mathcal{F}_{(A^T)^\xi \mathbf{x}_0}(t) = \mathcal{F}_{\mathbf{x}_0}(f)(t+\xi) .$$
    \end{proof}
\end{prop}

This proposition suggests that by grouping together terms in the expansion, we may obtain a decomposition of $f$ into smooth functions on $M$.

\begin{defn}
    Let $\orb_{\mathbf{y}}$ denote the orbit of the point $\mathbf y \in\mathbb{Z}^n$ being acted on by the group generated by the transpose matrix $A^T$. That is to say we have
    $$\orb_{\mathbf y} = \{(A^T)^\xi \mathbf y \,\mid\, \xi \in \mathbb{Z}\}. $$
    We use these orbits to partition $\mathbb{Z}^n$ and define $\mathcal{O}$ to be the set of all such orbits. 
\end{defn}
\begin{prop}\label{decomp}
    Any $f\in C^\infty(M)$ can be written as the series
    $$\sum_{\substack{\orb_{\mathbf y}\in\mathcal{O}\\ \abs{\orb_{\mathbf y}}=\infty}}\left(\sum_{\xi\in\mathbb{Z}} \mathcal{F}_{\mathbf y}(f)(t+\xi)e^{2\pi i \mathbf y \cdot A^\xi \mathbf x}\right) + \sum_{\substack{\orb_{\mathbf y}\in\mathcal{O}\\ \abs{\orb_{\mathbf y}}=N<\infty}}\left(\sum_{\xi = 0}^{N-1}\mathcal{F}_{\mathbf y}(f)(t+\xi)e^{2\pi i \mathbf y \cdot A^\xi \mathbf x}\right) $$
    and we have
    $$\left(\sum_{\xi\in\mathbb{Z}} \mathcal{F}_{\mathbf y}(f)(t+\xi)e^{2\pi i \mathbf y \cdot A^\xi \mathbf x}\right)\in C^\infty(M) $$
    $$\left(\sum_{\xi = 0}^{N-1}\mathcal{F}_{\mathbf y}(f)(t+\xi)e^{2\pi i \mathbf y \cdot A^\xi \mathbf x}\right)\in C^\infty(M) $$
    in the cases where $\mathbf{y} \in \mathbb{Z}^n$ satisfies $\abs{\orb_{\mathbf y}} =\infty$, respectively $\abs{\orb_{\mathbf y}} = N<\infty$.
    \begin{proof}
          By partitioning $\mathbb{Z}^n$ into the orbits $\orb_{\mathbf y}$ we can write
          $$\sum_{\mathbf{x}_0 \in \mathbb{Z}^n}\mathcal{F}_{\mathbf{x}_0}(f)(t)e^{2\pi i \mathbf{x}_0 \cdot \mathbf x} = \sum_{\orb_{\mathbf y} \in \mathcal{O}} \sum_{\mathbf{x}_0 \in \orb_{\mathbf y}}\mathcal{F}_{\mathbf{x}_0}(f)(t)e^{2\pi i \mathbf{x}_0 \cdot \mathbf x} .$$
          Then by Proposition \ref{function on M}, if we have $\mathbf{x}_0= (A^T)^\xi \mathbf y$ for some $\xi \in \mathbb{Z}$, then we can write
          $$\mathcal{F}_{\mathbf{x}_0}(f)(t) = \mathcal{F}_{\mathbf y}(f)(t+\xi) $$
          and thus
          $$\sum_{\mathbf{x}_0 \in \orb_{\mathbf y}}\mathcal{F}_{\mathbf{x}_0}(f)(t)e^{2\pi i \mathbf{x}_0 \cdot \mathbf x}= \sum_{\xi}\mathcal{F}_{\mathbf y}(f)(t+\xi)e^{2\pi i \mathbf y \cdot A^\xi \mathbf x} $$
          with $\xi$ ranging over different values depending on the size of $\orb_{\mathbf y}$.
    \end{proof}
\end{prop}

In the case when $\abs{\orb_{\mathbf y}}=N$ for some $N<\infty$ the function $\mathcal{F}_{\mathbf{y}}(f)$ is periodic with period $N$, and so we can further decompose it as follows

\begin{prop}\label{G decomp}
    Given $f\in C^\infty(M)$ and any $\mathbf y \in \mathbb{Z}^n$ such that $\abs{\orb_{\mathbf y}}=N<\infty$, we can write 
    $$\mathcal{F}_{\mathbf y}(f)(t) = \sum_{t_0 \in\mathbb{Z}}\mathcal{G}_{t_0,\mathbf y}(f) e^{\frac{2\pi i t_0 t}{N}} $$
    where $\mathcal{G}_{t_0, \mathbf y}\in\mathbb{C}$ is defined by
    $$\mathcal{G}_{t_0,\mathbf y}(f) = \frac 1N \int_0^N \mathcal{F}_{\mathbf y}(f)(t)e^{-\frac{2\pi i t_0 t}{N}}\, dt .$$
    \begin{proof}
          This is simply the Fourier expansion of the periodic function $\mathcal{F}_{\mathbf y}(f)(t)$.
    \end{proof}
\end{prop}
\begin{cor}
    In the decomposition of $f$ in Proposition \ref{decomp}, the summand
    $$ \left(\sum_{\xi = 0}^{N-1}\mathcal{F}_{\mathbf y}(f)(t+\xi)e^{2\pi i \mathbf y \cdot A^\xi \mathbf x}\right)\in C^\infty(M)$$
    can be further decomposed into
    $$\sum_{t_0\in \mathbb{Z}}\left(\mathcal{G}_{t_0, \mathbf y}(f)e^{2\pi i \frac{t_0 t}{N}}\sum_{\xi=0}^{N-1}e^{2\pi i \left(\frac {t_0 \xi}{N} + \mathbf y \cdot A^\xi \mathbf x\right) }\right) $$
    such that each term
    $$\left(\mathcal{G}_{t_0, \mathbf y}(f)e^{2\pi i \frac{t_0 t}{N}}\sum_{\xi=0}^{N-1}e^{2\pi i \left(\frac {t_0 \xi}{N} + \mathbf y \cdot A^\xi \mathbf x\right) }\right) $$
    is itself a smooth function on $M$.
    \begin{proof}
          This result is achieved by substituting the expression for $\mathcal{F}_{\mathbf y}(f)$ in Proposition \ref{G decomp} into the summand. That the terms of the decomposition are themselves smooth functions on $M$ can be verified through the use of Proposition \ref{function on M}.
    \end{proof}
\end{cor}

In the case when $\abs{\orb_{\mathbf y}}=\infty$ there does not seem to be any further useful decomposition of $\mathcal{F}$, however there are additional properties which $\mathcal{F}$ must satisfy.
\begin{prop}
    For any $f\in C^\infty(M)$ and any $\mathbf{y}\in \mathbb{Z}^n$ such that $\abs{\orb_{\mathbf y}}=\infty$, we require that all derivatives of $\mathcal{F}_{\mathbf y}(f)(t)$ tend to zero as $t\rightarrow \pm\infty$ faster than any power of $\abs{(A^T)^\xi \mathbf y}$ grows as $\xi \rightarrow \pm \infty$. Specifically, for any compact set $K\subset \mathbb{R}$ we require
    $$ \sup_{\substack{t\in K\\ \xi \in \mathbb{Z}}}\abs{\norm{(A^T)^\xi \mathbf y}^p \frac{d^q}{dt^q}\mathcal{F}_{\mathbf y}(f)(t+\xi)}<\infty $$
    for all $p,q\in \mathbb{N}$.
    \begin{proof}
        First, note that given any smooth function $f\in C^\infty(\mathbb{R}^{n+1})$ satisfying \eqref{periodicity}, all its derivatives must be bounded over any compact $\tilde{K}\subset \mathbb{R}^{n+1}$. If we take $\tilde{K}=[0,1]^n \times K$, we see that the Fourier coefficients $\mathcal{F}_{\mathbf{x}_0}$ of all the derivatives of $f$ must be bounded for $t$ ranging over $K$. Importantly, this bound is independent of $\mathbf{x}_0 \in \mathbb{Z}^n$.
        
        The Fourier coefficients of the derivatives of $f$ can take the form of $M(\mathbf{x}_0) \frac{d^q}{dt^q}\left(\mathcal{F}_{\mathbf{x}_0}(f)(t)\right)$ for any monomial $M$ and any $q\in\mathbb{N}$. This means for all monomials $M$ and all $q\in \mathbb{N}$ we require
        $$\sup_{\substack{t\in K\\\mathbf x_0 \in \mathbb{Z}^n}}\abs{M(\mathbf x_0) \frac{d^q}{dt^q}\mathcal{F}_{\mathbf x_0}(f)(t)}<\infty$$
        and thus if we restrict our attention to $\mathbf{x}_0 \in \orb_{\mathbf y}$ we require
        $$\sup_{\substack{t\in K\\\mathbf x_0 \in \orb_{\mathbf y}}}\abs{M(\mathbf x_0) \frac{d^q}{dt^q}\mathcal{F}_{\mathbf x_0}(f)(t)} = \sup_{\substack{t\in K\\ \xi \in \mathbb{Z}}}\abs{M((A^T)^\xi \mathbf y) \frac{d^q}{dt^q}\mathcal{F}_{\mathbf y}(f)(t+\xi)}<\infty.$$
        $M(\mathbf{x}_0)$ can then be chosen to be $\norm{\mathbf{x}_0}^p$ for arbitrarily large $p\in \mathbb{N}$, giving us the desired result.
    \end{proof}
\end{prop}
\begin{cor} \label{Schwartz}
    For any $f\in C^\infty(M)$ and any $\mathbf{y}\in \mathbb{Z}^n$ such that $\abs{\orb_{\mathbf y}}=\infty$, we require that $\mathcal{F}_{\mathbf y}(f)(t) \in \mathcal{S}(\mathbb{R})$. Here $\mathcal{S}(\mathbb{R})$ denotes the space of Schwartz functions
    $$\mathcal{S}(\mathbb{R}) = \left\{h(t)\in C^\infty(\mathbb{R}) \,\middle| \, \sup_{t\in\mathbb{R}} \abs{t^p \frac{d^q}{dt^q}h(t)}<\infty, \text{ for all } p,q \in \mathbb{N}  \right\}. $$
    \begin{proof}
          If $\abs{\orb_{\mathbf y}}=\infty$ then $ \norm{(A^T)^\xi \mathbf y}$ must blow up as $\xi \rightarrow \pm \infty$, since an infinite orbit cannot repeat the same point twice. Furthermore, since the number of lattice points within a bounded region of $\mathbb{Z}^n$ grows like $R^2$ with the radius $R$ of the region, it must be the case that $ \norm{(A^T)^\xi \mathbf y}$ blows up at least as fast as $\abs{\xi}^{\frac 12}$. Substituting this speed of growth into the above proposition gives the definition of $\mathcal{S}(\mathbb{R})$. 
          
          Note that if $ \norm{(A^T)^\xi \mathbf y}$ blows up faster than polynomially, then the Proposition yields an even stricter condition on $\mathcal{F}_{\mathbf y}$ than Schwartz.   
    \end{proof}
\end{cor}

\begin{thm}\label{L2 decomp}
  The space of $L^2$ functions on $M$ decomposes in the following way.
  $$L^2(M) = \left(\widehat{\bigoplus_{\substack{\orb_{\mathbf y}
  \in \mathcal{O}\\ \abs{\orb_{\mathbf y}}=\infty}}}\mathcal{H}_{\mathbf y}\right) \oplus  \left(\widehat{\bigoplus_{\substack{\orb_{\mathbf y}
  \in \mathcal{O}\\ \abs{\orb_{\mathbf y}}= N<\infty}}}\widehat{\bigoplus_{t_0 \in \mathbb{Z}}}\mathcal{H}_{t_0, \mathbf y}\right), $$
  where
  $$\mathcal{H}_{\mathbf y} = \left\{ \sum_{\xi\in\mathbb{Z}} f(t+\xi)e^{2\pi i \mathbf y \cdot A^\xi \mathbf x} \,\middle|\, f \in L^2(\mathbb{R})  \right\}  $$
  and
  $$\mathcal{H}_{t_0, \mathbf y} = \left\{ C e^{2\pi i \frac{t_0 t}{N}}\sum_{\xi = 0}^{N-1} e^{2\pi i\left(\frac{t_0 \xi}{N}+\mathbf y \cdot A^{\xi} \mathbf x\right)}\, \middle|\, C\in \mathbb{C} \,\right\}. $$
  Here $\hat \oplus$ denotes the direct sum followed by the closure with respect to the $L^2$ norm.
  \begin{proof}
        From Propositions \ref{decomp} and \ref{G decomp} we see that any smooth function can be decomposed in the way described above and so, since $L^2(M)$ is the closure of $C^{\infty}(M)$ with respect to the $L^2$ norm, we obtain the desired result.  
  \end{proof}
\end{thm}

It will be useful now to consider what the orbits of $\mathbf y \in \mathbb{Z}^n$ actually look like. In particular, when exactly is $\abs{\orb_{\mathbf{y}}}<\infty$. First, we define the generalised eigenvectors of $A$.

\begin{defn}
    Let $\lambda_1, \dots, \lambda_k\in\mathbb{C}$ be the eigenvalues of $A \in GL_n(\mathbb{Z})$ with values repeated for geometric multiplicity. Then any $n$ linearly independent vectors $\mathbf{v}_{i,j}\in\mathbb{C}^n$ with $i=1,\dots, k$ and $j=1,\dots, m_i$ such that
        $$(A-\lambda_i)^j \mathbf{v}_{i,j} = 0 \quad \text{ but }\quad (A-\lambda_i)^{j-1} \mathbf{v}_{i,j} \neq 0 $$
    are called generalised eigenvectors of $A$. Note that when $j=1$ we just have the standard eigenvectors of $A$. Furthermore, we can make a choice of $\mathbf{v}_{i,j}$ so that when $i$ is fixed, the sequence $\mathbf{v}_{i,1}, \mathbf{v}_{i,2}, \dots , \mathbf{v}_{i,m_i}$ forms a Jordan chain of length $m_i$. This means for all $j\neq 1$ we have
        \begin{equation}
            (A-\lambda_i)\mathbf{v}_{i,j}=\mathbf{v}_{i,j-1} 
        \label{jordan1}
        \end{equation}
        and for $j=1$ we have
         \begin{equation}
            (A-\lambda_i)\mathbf{v}_{i,1}=0.
        \label{jordan2}
        \end{equation}

\end{defn}

These $\mathbf{v}_{i,j}$ can be used to describe when the orbit of the group generated by $A^T$ acting on $\mathbf y \in \mathbb{Z}^n$ is finite.

\begin{prop}
    Let $\mathbf{v}_{i,j}$ be the generalised eigenvectors of $A\in GL_{n}(\mathbb{Z})$ as defined above, with corresponding eigenvalues $\lambda_1 , \dots , \lambda_k$.  Given $\mathbf{y}\in \mathbb{Z}^n$, if $\abs{\orb_{\mathbf{y}}}=N<\infty$ it must be the case that $\mathbf{v}_{i,j} \cdot \mathbf{y} =0$ except for when $i$ and $j$ are chosen such that $\lambda_i^N = 1$ and $j=m_i$
    \begin{proof}
        If $\orb_{\mathbf y}$ is a finite subset of $\mathbb{Z}^n$, then $(A^T)^\xi \mathbf y$ 
        must be bounded over $\xi \in \mathbb{Z}$. 
        This means $\mathbf{v}_{i,j}\cdot ((A^T)^\xi \mathbf y)= (A^\xi \mathbf{v}_{i,j})\cdot \mathbf{y}$ must be bounded over $\xi \in \mathbb{Z}$ for all $\mathbf{v}_{i,j}$.
        
        From \eqref{jordan2} we know that
        $A \mathbf{v}_{i,1} = \lambda \mathbf{v}_{i,1}$ and thus
        $$A^\xi \mathbf{v}_{i,1} \cdot \mathbf y = \lambda^\xi \mathbf{v}_{i,1}\cdot  \mathbf y.$$
        But if $\abs{\lambda_i}>1$ then $\lambda_i^\xi$ will blow up as $\xi \rightarrow \infty$ and if $\abs{\lambda_i}<1$ then it will blow up as $\xi \rightarrow -\infty$. From this we conclude that $\abs{\orb_{\mathbf y}}<\infty$ only if $\mathbf{v}_{i,1}\cdot  \mathbf y = 0$ for all $i$ such that $\abs{\lambda_i}\neq 1$. Rewriting $\eqref{jordan1}$ as $A \mathbf{v}_{i,j} = \lambda \mathbf{v}_{i,j}+\mathbf{v}_{i,j-1}$ and using $\mathbf{v}_{i,1}\cdot \mathbf y =0$ we can apply the above argument again to prove the same result for $\mathbf{v}_{i,2}$. In fact, continuing by induction, we see that $\abs{\orb_{\mathbf y}}$ is finite only if $\mathbf{v}_{i,j}\cdot  \mathbf y = 0$ for all $i$ and $j$ such that $\abs{\lambda_i}\neq 0$.
        
        Now, consider the case when $\abs{\lambda_i}=1$. From \eqref{jordan1} we can see that when $m_i\geq 2$ then 
        $$A^\xi \mathbf{v}_{i,2} = \lambda^\xi \mathbf{v}_{i,2}+ \xi \lambda^{\xi-1}\mathbf{v}_{i,1} .$$
        This means $A^\xi \mathbf{v}_{i,2} \cdot \mathbf{y}$ will blow up as $\xi \rightarrow \pm\infty$ unless $\mathbf{v}_{i,1}\cdot \mathbf{y} = 0$. Similarly, if $\mathbf{v}_{i,1}\cdot \mathbf{y} = 0$ then the same argument works to show $A^\xi \mathbf{v}_{i,3} \cdot \mathbf{y} $ will blow up unless $\mathbf{v}_{i,2}\cdot \mathbf{y} = 0$, provided $m_i \geq 3$. Repeating this procedure, we find that $\abs{\orb_{\mathbf y}}<\infty$ implies that $\mathbf{v}_{i,j}\cdot  \mathbf y = 0$ for all $i$ and $j$ such that $\abs{\lambda_i}=1 $ and $j < m_i$
        
        Finally, it remains to consider the case of $\mathbf{v}_{i,m_i}$. If $\abs{\orb_{\mathbf y}}=N$ then we know that $(A^T)^N \mathbf{y} = \mathbf y$, and also we have shown that $\mathbf{v}_{i,j}\cdot  \mathbf y = 0$ for all $j \neq m_i$. The following must therefore hold.
        \begin{align*}
            \mathbf{v}_{i,m_i} \cdot \mathbf y &= \mathbf{v}_{i,m_i} \cdot (A^T)^N\mathbf y\\ 
            &= A^N \mathbf{v}_{i,m_i} \cdot \mathbf y\\ &= \lambda_i^N \mathbf{v}_{i,m_i} \cdot \mathbf y
        \end{align*}
        Thus $\abs{\orb_{\mathbf y}}=N$ requires that for all $i$, either $\mathbf{v}_{i,m_i} \cdot \mathbf y =0$ or $\lambda_i^N = 1$
    \end{proof}
\end{prop}

\begin{cor} \label{theta}
    Whenever $\abs{\orb_{\mathbf y}}=N<\infty$, it holds that
    $$ A \mathbf{v}_{i,j}\cdot \mathbf y =  \begin{cases}
        e^{2\pi i \theta_i}\mathbf{v}_{i,j}\cdot \mathbf y & \text{ if } \lambda_i^N = 1 \text{ and } j=m_i\\
        0 & \text{ otherwise }
    \end{cases}
    $$
    where $\theta_i \in\mathbb{Q}\cap (-\frac 12 , \frac 12 ] $ is some rational number depending on $i$ satisfying $N\theta_i \in \mathbb{Z}$. 
\end{cor}

\subsection{Properties of the decomposition}
We would now like to consider some of the properties of this decomposition, which will be useful when considering the example in the following section. But in order to do this we must first construct a special frame on $M$.

\begin{defn} \label{matrix power}
    Given any invertible matrix $A\in GL_n(\mathbb{Z})$, then for some choice of matrix logarithm $\ln A$ we can define the power $A^t := e^{t \ln A}$
    for all $t \in \mathbb{R}$. Note that such a logarithm always exists, but may be complex valued.
    
    Throughout this paper, the choice of $\ln A$ will always be made such that 
        $$ A^t \mathbf{v}_{i,j}\cdot \mathbf y =  \begin{cases}
        e^{2\pi i \theta_i t}\mathbf{v}_{i,j}\cdot \mathbf y & \text{ if } \lambda_i^N = 1 \text{ and } j=m_i\\
        0 & \text{ otherwise }
    \end{cases}
    $$
    for $\theta_i \in \mathbb{Q}\cap (-\frac 12,\frac 12 ]$
\end{defn}

Using the generalised eigenvectors of $A$ given by $\mathbf{v}_{i,j}$, a smooth frame for the complexified tangent bundle of $M$ can be given by
$$\epsilon_0 = \frac{\partial}{\partial t} \quad \quad \epsilon_{i,j}= A^t \mathbf{v}_{i,j} \cdot \nabla_{\mathbf{x}}. $$
Here we are using $\nabla_{\mathbf{x}}=(\frac{\partial}{\partial x_1},\frac{\partial}{\partial x_2},\dots, \frac{\partial}{\partial x_n})$ to denote the gradient excluding the variable $t$. We verify that this is indeed a well defined frame on $M$ in following proposition.

\begin{prop}
    Viewing $M$ as a torus bundle over $S^1$, any smooth frame of the complexified tangent bundle on a single fibre may be extended to a smooth frame on all of $M$.
    
    \begin{proof}
        We can assume, without loss of generality, that we are starting with a frame on the $t=0$ fibre, where $t$ is parametrising the base space $S^1$, as in the definition of $M$ \eqref{identification}. 
        
        Let $\mathbf{a_1},\mathbf{a_2}, \dots , \mathbf{a_n}:\mathbb{R}^n / \mathbb{Z}^n \rightarrow \mathbb{C}^n$ be smooth maps sending each point $\mathbf{x}\in \mathbb{T}^n$ to $n$ linearly independent vectors. Then the collection $\{\mathbf{a}_i \cdot \nabla_{\mathbf x} \}_{i=1,\dots,n }$ defines a general frame on the $t=0$ fibre.
        A frame for $T_\mathbb{C}M$ is then given by $u_0 = \frac{\partial}{\partial t}$ and  $u_i=A^t \mathbf{a}_i \cdot \nabla_{\mathbf x}$ with $i=1,\dots,n$
        
        These are indeed all well defined vector fields on $M$, in particular they do not conflict with the identification of points given in $\eqref{identification}$. To check the first identification, simply note that the maps $\mathbf{a}_i (\mathbf{x})$ are defined on the torus. For the second we consider the map
        $$\phi_{\xi}:\begin{pmatrix}
        t\\ \mathbf x 
        \end{pmatrix}\mapsto \begin{pmatrix}
        t+\xi \\ A^\xi \mathbf x
        \end{pmatrix} $$
        with $\xi \in \mathbb{Z}$ and try to show that $u_i$ are invariant under the pushforward. Certainly this is true of $\frac{\partial}{\partial t}$, and we also know that, for $i = 1,\dots, n$, we have
        \begin{align*}
            (\phi_\xi)_* (\mathbf{e}_i \cdot \nabla_{\mathbf{x}}) &=  (\phi_\xi)_*  \frac{\partial}{\partial x_i}\\
            &= A^\xi \mathbf{e}_i \cdot \nabla_{\mathbf x}
        \end{align*}  
        with $\mathbf{e}_i$ signifying the standard basis vector $(0,\dots , 1,\dots , 0)$ with a 1 in the $i^{th}$ position. Therefore
        \begin{align*}
            (\phi_\xi )_* u_i(t) &= (\phi_\xi )_* (A^t \mathbf{a}_i \cdot \nabla_{\mathbf x}) \\
            &= A^{t+\xi} \mathbf{a}_i \cdot \nabla_{\mathbf x}\\
            &= u_i(t+\xi).
        \end{align*} 
        
    \end{proof}
\end{prop}

    It should be noted that if $A$ has a real-valued logarithm and we choose $\mathbf{a}_{i}$ to be maps into $\mathbb{R}^n$, then the construction in the above proof will give us a smooth frame on the standard, non-complexified tangent bundle.
    
\begin{prop}\label{F}
    Given any $\mathbf y \in \mathbb{Z}^n$ and any $f\in C^\infty(M)$, $\mathcal{F}_{\mathbf y}$ has the properties
    \enumerate[label=\roman*)]
        \item $$\mathcal{F}_{\mathbf y}(\epsilon_0 f)(t)=\epsilon_0 \,\mathcal{F}_{\mathbf y}(f)(t), $$
        \item $$\mathcal{F}_{\mathbf y}(\epsilon_{i,j} f)(t)= 2\pi i A^t \mathbf{v}_{i,j} \cdot \mathbf y \,\mathcal{F}_{\mathbf y}(f).  $$
    \begin{proof}
         Since $\mathcal{F}_{\mathbf y}(f)$ is just one of the Fourier coefficients of $f$ in the standard expansion
         $$f(t,\mathbf x) = \sum_{\mathbf y \in\mathbb{Z}^n} \mathcal{F}_{\mathbf y}(f)(t) e^{2\pi i \mathbf y \cdot \mathbf x} $$
         this proposition is simply restating results from classical Fourier analysis,
    \end{proof}
\end{prop}

\begin{prop}\label{G}
    Given any $\mathbf{y} \in \mathbb{Z}^n$ such that $\abs{\orb_{\mathbf y}}=N<\infty$ and any $f\in C^\infty(M), \mathcal{G}_{t_0, \mathbf y}$ has the properties
    \enumerate[label=\roman*)]
        \item $$\mathcal{G}_{t_0, \mathbf y}(\epsilon_0 f)(t)=2\pi i \frac{t_0}{N} \,\mathcal{G}_{t_0, \mathbf y}(f), $$
        \item $$\mathcal{G}_{t_0, \mathbf y}(\epsilon_{i,j} f)= \begin{cases}2\pi i \mathbf{v}_{i,j} \cdot \mathbf y \, \mathcal{G}_{t_0 + N\theta_i ,\mathbf y}(f) & \text{ if } \lambda_i^N = 1 \text{ and } j=m_i\\
        0 & \text{ otherwise }
        \end{cases}.  $$
    With $\theta_i$ defined as in Corollary \ref{theta}.
    
    \begin{proof}
          For part $i)$, we make use of the result $i)$ in the previous proposition along with the definition of $\mathcal{G}_{t_0, \mathbf y}$ to write
          \begin{align*}
              \mathcal{G}_{t_0, \mathbf y}(\epsilon_0 f)(t) &= \frac 1N \int_0^N \mathcal{F}_{\mathbf y}(\epsilon_0 f) e^{-\frac{2\pi i t_0 t}{N}}\,dt\\
              &= \frac 1N \int_0^N \left(\epsilon_0 \mathcal{F}_{\mathbf y}(f)\right) e^{-\frac{2\pi i t_0 t}{N}}\, dt.
          \end{align*}
          Then, since $\mathcal{F}_{\mathbf y}(f)(t)$ is periodic with period $N$, we can make use of integration by parts to get
          \begin{align*}
              \frac 1N \int_0^N \left(\epsilon_0 \mathcal{F}_{\mathbf y}(f)\right) e^{-\frac{2\pi i t_0 t}{N}}\, dt&=
              -\frac 1N \int_0^N \mathcal{F}_{\mathbf y}(f) \left(\epsilon_0 e^{-\frac{2\pi i t_0 t}{N}}\right)\, dt\\
              &= 2\pi i \frac{t_0}{N} \frac 1N \int_0^N \mathcal{F}_{\mathbf y}(f) e^{-\frac{2\pi i t_0 t}{N}}\, dt\\
              &= 2\pi i \frac{t_0}{N} \mathcal{G}_{t_0,\mathbf y}(f).
          \end{align*}
          
          For part $ii)$, we make use of the result $ii)$ in the previous proposition to write
          \begin{align*}
              \mathcal{G}_{t_0,\mathbf y}(\epsilon_{i,j}f) &= \frac 1N \int_0^N \mathcal{F}_{\mathbf y}(\epsilon_{i,j} f) e^{- \frac{2\pi i t_0 t}{N}}\, dt\\
              &= \frac 1N \int_0^N 2\pi i A^t \mathbf{v}_{i,j} \cdot \mathbf y \mathcal{F}_{\mathbf y}(f) e^{- \frac{2\pi i t_0 t}{N}}\, dt.
          \end{align*}
          Then, because of the way $A^t$ was defined in Definition \ref{matrix power}, we get
          $$
               \frac 1N \int_0^N 2\pi i A^t \mathbf{v}_{i,j} \cdot \mathbf y \mathcal{F}_{\mathbf y}(f) e^{- \frac{2\pi i t_0 t}{N}}\, dt =0   
         $$
         unless $\lambda_i^N = 1$ and $j=m_i$, in which case
         \begin{align*}
              \frac 1N \int_0^N 2\pi i A^t \mathbf{v}_{i,j} \cdot \mathbf y \mathcal{F}_{\mathbf y}(f) e^{- \frac{2\pi i t_0 t}{N}}\, dt
              &= 2\pi i \mathbf{v}_{i,j} \cdot \mathbf y  \frac 1N \int_0^N   \mathcal{F}_{\mathbf y}(f) e^{- \frac{2\pi i (t_0 +N\theta_i)t}{N}}\, dt\\
              &= 2\pi i\mathbf{v}_{i,j} \cdot \mathbf y \, \mathcal{G}_{t_0+N\theta_i, \mathbf y}(f).
         \end{align*}
    \end{proof}
\end{prop}

\section{Calculating $h^{0,1}$ using Harmonic Analysis}
It should be noted the Kodaira-Thurston manifold $KT^4$ can be viewed as a torus bundle over $S^1$ with $A= \begin{pmatrix}
1 & 0 & 0\\
0 & 1 & 1\\
0 & 0 & 1
\end{pmatrix} $. The calculations done on $KT^4$ in Subsection \ref{example KT4} can be thought of as an application of the above results, with the case when $n\neq 0$ corresponding to an infinite orbit and the case when $n = 0$ corresponding to a finite orbit of length 1.

In this section we will see what it looks like to use our decomposition to perform calculations on a manifold for which we have a finite orbit with length greater that 1. 
To that end, we will define $M$ by setting
$$A = \begin{pmatrix}
0 & 0 & 1\\
1 & 0 & 0\\
0 & 1 & 0
\end{pmatrix} $$
and identifying points in $\mathbb R^4$ by
$$\begin{pmatrix} t \\ \mathbf x \end{pmatrix}
\sim \begin{pmatrix} t \\ \mathbf x + \mathbf x_0 \end{pmatrix} \quad \text{ and } \quad 
\begin{pmatrix} t \\ \mathbf x \end{pmatrix} \sim 
\begin{pmatrix} t+\xi \\ A^\xi \mathbf x\end{pmatrix}
$$
for all $\mathbf x_0 \in \mathbb Z^3$ and all $\xi \in \mathbb Z$. 
The matrix $A$ has eigenvalues of 1, $e^{-\frac 23 \pi i}$ and $e^{\frac 23 \pi i}$ corresponding to eigenvectors $\begin{pmatrix} 1 \\ 1 \\ 1 \end{pmatrix}$,$\begin{pmatrix} e^{\frac 23 \pi i} \\ e^{-\frac 23 \pi i} \\ 1 \end{pmatrix}$ and $\begin{pmatrix} e^{\frac 23 \pi i} \\ e^{-\frac 23 \pi i} \\ 1 \end{pmatrix}$. We therefore define a smooth frame on the complexified tangent bundle by
$$\epsilon_0 = \der t  \quad \quad \epsilon_1 = \begin{pmatrix} 1 \\ 1 \\ 1 \end{pmatrix} \cdot \nabla_\mathbf x \quad \quad 
\epsilon_2 = e^{-\frac 23 \pi i t}\begin{pmatrix} e^{\frac 23 \pi i} \\ e^{-\frac 23 \pi i} \\ 1 \end{pmatrix} \cdot \nabla_\mathbf x
\quad \quad 
\epsilon_3 = e^{\frac 23 \pi i t}\begin{pmatrix} e^{-\frac 23 \pi i} \\ e^{\frac 23 \pi i} \\ 1 \end{pmatrix} \cdot \nabla_\mathbf x$$
where we define $\nabla_{\mathbf x} := \begin{pmatrix}
\partial_x\\ \partial_y\\ \partial_z
\end{pmatrix}$. The dual frame is given by
$$  
\epsilon^0 = \dif t, \quad \quad 
\epsilon^1 = \frac 13 \left(\dif x + \dif y + \dif z \right),$$
$$\epsilon^2 = \frac{e^{\frac 23 \pi i t}}{3}\left(e^{-\frac 23 \pi i }\dif x + e^{\frac 23 \pi i } \dif y + \dif z\right), \quad \quad
\epsilon^3 = \frac{e^{-\frac 23 \pi i t}}{3}\left(e^{\frac 23 \pi i }\dif x + e^{-\frac 23 \pi i } \dif y + \dif z\right).$$

Let an almost complex structure $J$ be defined by the mapping
$$\epsilon_0 \mapsto \frac 12 (\epsilon_2 + \epsilon_3) \quad \quad \text{ and } \quad \quad \epsilon_1 \mapsto -\frac {i}{2} (\epsilon_2 - \epsilon_3). $$
We can then find a pair of vectors fields spanning $T^{1,0}_pM$ at all points $p \in M$
$$V_1 = \frac 12 \left(\epsilon_0 - \frac i2 (\epsilon_2 + \epsilon_3) \right) \quad \quad V_2 = \frac 12 \left(\epsilon_1 + \frac 12 (\epsilon_2 - \epsilon_3)  \right)  $$
with dual $(1,0)$-forms given by
$$\phi^1 = \epsilon^0 + i(\epsilon^2 + \epsilon^3) \quad \quad \phi^2 = \epsilon^1 -\epsilon^2 + \epsilon^3  $$
which satisfy the structure equations
$$\dif \phi^1 = \frac {\pi}{6} \left(\phi^{12} - \phi^{1 \bar 2} - \phi^{2 \bar 1} - \phi^{\bar 1 \bar 2}\right) $$
$$\dif \phi^2 = \frac{\pi}{3} \phi^{1 \bar 1}.$$
The metric can be chosen so that $V_1$ and $V_2$ form a unitary basis. 

\subsection{Deriving the equations}
Let a general $(0,1)$-form be written as $s = f \bar\phi^1 + g \bar\phi^2$. The two requirements 
$$\bar \partial s = 0 \quad \quad \text{ and } \quad \quad \partial * s = 0 $$
which are equivalent to $s$ being $\bar\partial$-harmonic, give rise to the two PDEs
$$-(\bar V_2 - \frac \pi 6) f + \bar V_1 g = 0 $$
$$ V_1 f + V_2 g = 0 .$$
We will now try taking a Fourier expansion. If $f$ is a smooth function on $M$ we know from Proposition \ref{function on M} that we can write it as
$$ f(t,x,y,z) = \sum_{l,m,n\in \mathbb Z} \mathcal{F}_{l,m,n}(f)(t)e^{2\pi i (lx+my+nz)}.$$
with $\mathcal{F}_{l,m,n}(f)$ satisfying the property
$$\mathcal{F}_{m,n,l} (t) = \mathcal{F}_{l,m,n}(f)(t+\xi) $$
for all $\xi \in \mathbb Z$.
This gives us two cases: if $l=m=n$ then  $\mathcal{F}_{l,m,n}(f)$ is periodic with period length 1, and otherwise it is periodic with period length 3. So in both cases we can further expand the function with respect to the variable $t$.
\subsubsection{The $l=m=n$ case}
Here we have a standard Fourier expansion in all 4 variables, with
$$\mathcal{F}_{n,n,n}(f) = \sum_{k\in \mathbb Z}\mathcal{G}_{k,n,n,n}(f)(t)e^{2\pi i kt} .$$
In this case by Proposition \ref{G} we see that $\mathcal{G}_{k,n,n,n}$ satisfies the properties
\begin{align*}
    \mathcal{G}_{k,n,n,n}(\epsilon_0 f) &= 2\pi i k \mathcal{G}_{k,n,n,n}(f),\\ 
    \mathcal{G}_{k,n,n,n}(\epsilon_1 f) &= 2\pi i (n+n+n) \mathcal{G}_{k,n,n,n}(f) = 6\pi i n \mathcal{G}_{k,n,n,n}(f) ,\\ 
    \mathcal{G}_{k,n,n,n}(\epsilon_2 f) &= 2\pi i (e^{\frac 23 \pi i}n + e^{-\frac 23 \pi i}n + n)\mathcal{G}_{k,n,n,n}(f) = 0\\ \mathcal{G}_{k,n,n,n}(\epsilon_3 f) &= 2\pi i (e^{-\frac 23 \pi i}n + e^{\frac 23 \pi i}n + n)\mathcal{G}_{k,n,n,n}(f) = 0
\end{align*}
which we can use to rewrite our two PDEs into the form
$$\begin{pmatrix}
k & 3n\\
-3n-\frac i6 & k
\end{pmatrix} 
\begin{pmatrix}
\mathcal{G}(f)\\
\mathcal{G}(g)
\end{pmatrix} = 0.$$
This has non-trivial solutions if and only if the matrix has zero determinant, \textit{i.e.}
$$k^2 + 9n^2 + 3n\frac i6 = 0 $$
which is only the case when $k=n=0$. Corresponding to this case we have the solution
$$f = 0 \quad \quad g = const. $$

\subsubsection{The second case}
In this case $\mathcal{F}_{l,m,n}(f)$ is still periodic, but with period 3 and so our expansion now looks like 
$$\mathcal{F}_{l,m,n}(f)(t) = \sum_{k \in\mathbb{Z}}\mathcal{G}_{k,l,m,n}(f) e^{\frac{2\pi i k t}{3}}. $$
For the sake of notational simplicity we define
$$\alpha_1 = l+m+n \quad \quad \alpha_2 = e^{\frac 23 \pi i}l + e^{-\frac 23 \pi i}m + n \quad \quad \alpha_3 = e^{-\frac 23 \pi i}l + e^{\frac 23 \pi i}m + n $$
and also we will use $\mathcal{G}_{k}(f)$ to denote $\mathcal{G}_{k,l,m,n}(f)$. Then by Proposition \ref{G} we can say that $\mathcal{G}_{k}(f)$ satisfies
\begin{align*}
    \mathcal{G}_{k}(\epsilon_0 f) &= 2\pi i k \mathcal{G}_{k}(f),\\ 
    \mathcal{G}_{k}(\epsilon_1 f) &= 2\pi i \alpha_1 \mathcal{G}_{k}(f),\\ 
    \mathcal{G}_{k}(\epsilon_2 f) &= 2\pi i \alpha_2\mathcal{G}_{k-1}(f)\\ \mathcal{G}_{k}(\epsilon_3 f) &= 2\pi i \alpha_3\mathcal{G}_{k+1}(f).
\end{align*}
Applying these properties to our two PDEs gives us a pair of equations
\begin{align*}
   & \frac{\alpha_3}{2}\mathcal{G}_{k-1}(f) -(\alpha_1+\frac i6)\mathcal{G}_{k}(f) - \frac{\alpha_2}{2}\mathcal{G}_{k+1}(f)\\
    &+\frac{i\alpha_3}{2}\mathcal{G}_{k-1}(g)+\frac k3 \mathcal{G}_k (g) +\frac{i\alpha_2}{2}\mathcal{G}_{k+1}(g)=0\\
\end{align*}
\begin{align*}
    -&\frac{i\alpha_3}{2}\mathcal{G}_{k-1}(f)+ \frac k3 \mathcal{G}_k (f) -\frac{i\alpha_2}{2}\mathcal{G}_{k+1}(f)\\
    &+ \frac{\alpha_3}{2}\mathcal{G}_{k-1}(g) +\alpha_1 \mathcal{G}_k (g) -\frac{\alpha_2}{2}\mathcal{G}_{k+1}(g)=0.
\end{align*}
By choosing to cancel either the terms $\mathcal{G}_{k-1}(f)$ \& $\mathcal{G}_{k-1}(g)$ or the terms $\mathcal{G}_{k+1}(f)$ \& $\mathcal{G}_{k+1}(g)$ we can simplify to the pair of equations
$$\left(\frac k3 + \frac 16 - i \alpha_1 \right) \mathcal{G}_k(f)- i\alpha_2 \mathcal{G}_{k+1}(f) + i\left(\frac k3 -i\alpha_1\right)\mathcal{G}_k(g)+\alpha_2 \mathcal{G}_{k+1}(g)=0   $$
$$-i\alpha_3 \mathcal{G}_{k-1}(f) + \left(\frac k3 -\frac 16 + i\alpha_1\right)\mathcal{G}_k(f) + \alpha_3 \mathcal{G}_{k-1}(g) -i\left(\frac k3 + i\alpha_1\right)\mathcal{G}_k(g) = 0.  $$
Evaluating the second of these at $k+1$ instead of $k$ we can cancel either the $\mathcal{G}_{k+1}(f)$ term or the $\mathcal{G}_{k+1}(g)$ term. In this way we can write our equations as the recurrence relation
$$
\begin{pmatrix}
\mathcal{G}_{k+1}(f)\\ \mathcal{G}_{k+1}(g)
\end{pmatrix}
= 
\frac{6}{\left(4k + 3 + 12i\alpha_1\right)} B_k
\begin{pmatrix}
\mathcal{G}_k(f)\\ \mathcal{G}_k(g)
\end{pmatrix}$$
where 
$$B_k = \begin{pmatrix}
    -i\left[\left(\frac k3 + \frac 16\right)\left(\frac k3 + \frac 13\right)+\alpha_1^2 -\frac 16 i\alpha_1 - \alpha_2 \alpha_3 \right]&
    \left[\frac k3 \left(\frac k3 + \frac 13\right)+\alpha_1^2 -\frac 13 i\alpha_1 - \alpha_2 \alpha_3 \right]\\
    -\left[\left(\frac k3 + \frac 16\right)^2+\alpha_1^2 + \alpha_2 \alpha_3 \right] &
    -i\left[\frac k3\left(\frac k3 + \frac 16\right)+\alpha_1^2 -\frac 16 i\alpha_1 + \alpha_2 \alpha_3 \right]
\end{pmatrix}, $$
and so the values of $\mathcal{G}_{k}(f)$ and $\mathcal{G}_{k}(g)$ for all $k\in \mathbb{Z}$ are determined by a choice for $\mathcal{G}_{0}(f)$ and $\mathcal{G}_{0}(g)$.
Since we are looking for smooth solutions $f$ and $g$, we require that $\mathcal{F}_{l,m,n}(f)(t) = \sum_{k \in\mathbb{Z}}\mathcal{G}_{k,l,m,n}(f) e^{\frac{2\pi i k t}{3}} $ be smooth, and likewise for $\mathcal{F}_{l,m,n}(g)(t)$. This is equivalent to asking that the sequences $\mathcal{G}_{k}(f)$ and $\mathcal{G}_{k}(g)$ are Schwartz, \textit{i.e.} they are contained in 
$$\mathcal{S}(\mathbb{Z}) = \left\{(a_k)_{k \in \mathbb{Z}}  \,\middle|\, \sup_{k \in \mathbb{Z}}\abs{k^p a_k}<\infty \, \textit{ for all } p \in \mathbb{N}\right\}. $$

\end{document}